\documentclass[11pt]{article}

\usepackage{amsfonts}
\usepackage{amssymb}
\usepackage{mathrsfs}
\usepackage{amsmath}
\usepackage{amscd}
\usepackage{epsfig}
\usepackage{graphicx,dblfloatfix,caption}
\usepackage{listings}

\newtheorem{theorem}{Theorem}
\newtheorem{definition}{Definition}
\newtheorem{lemma}{Lemma}
\newtheorem{conjecture}{Conjecture}

\numberwithin{equation}{section}
\numberwithin{definition}{section}
\numberwithin{lemma}{section}

\begin{document}

\title{Rigidity of Infinite Hexagonal Triangulation of the Plane}

\author{
Tianqi Wu
\thanks{Mathematical Sciences Center, Tsinghua University, Beijing 100084 China. Email: mike890505@gmail.com.}
\and 
Xianfeng Gu
\thanks{Department of Computer Science, Stony Brook University, New York 11794 USA. Email: gu@cs.stonybrook.edu.}
\and
Jian Sun
\thanks{Mathematical Sciences Center, Tsinghua University, Beijing 100084 China. Email: jsun@math.tsinghua.edu.cn.}
}

\date{}
\maketitle

\begin{abstract}
In the paper, we consider the rigidity problem of the infinite
hexagonal triangulation of the plane under the piecewise linear
conformal changes introduced  by Luo in \cite{luo_yamabe}.  
Our result shows that if a geometric hexagonal triangulation of the plane is
PL conformal to the  regular hexagonal triangulation and all inner
angles are in $[\delta, \pi/2 -\delta]$ for any constant $\delta > 0$, 
then it is the regular hexagonal triangulation. This partially solves 
a conjecture of Luo~\cite{luo1}. The proof uses the concept of
\emph{quasi-harmonic} functions to unfold the properties of the
mesh.
\end{abstract}



\newpage

\section{Introduction}
\subsection{PL conformal}
Given a smooth manifold $M$, two Riemannian metrics $g$ and $\tilde{g}$ are called conformally equivalent if
\begin{equation}
\tilde{g} = e^{2\lambda} g
\end{equation}
where $\lambda \in C^\infty(M)$ is called a conformal factor. In
the discrete setting,  Luo~\cite{luo_yamabe} introduced a
notion of  PL conformal equivalence of two piecewise linear
polyhedral metrics in any dimension and developed a variational
principle for PL conformality for triangulated surfaces with PL
metrics ~\cite{luo_yamabe}. Specifically, Suppose the surface
$\Sigma$ has a triangulation $T$, i.e., a CW complex whose faces
are triangles which are glued edge-to-edge by isometrics. We
denote such triangulated surface by $(\Sigma, T)$,  and its sets
of vertices, edges, and triangles of $T$ by $V$, $E$ and $F$,
respectively. Two triangulated surfaces $(\Sigma, T)$ and
$(\tilde{\Sigma}, \tilde{T})$ are called combinatorially
equivalent if there is an homeomorophism between $\Sigma$ and
$\tilde{\Sigma}$ preserving the triangulation. For simplicity, we
use the same notation to denote two combinatorially equivalent
triangulated surfaces when the homeomorphism is not relevant or
clear. Recall that an \it piecewise linear metric \rm (or simply
PL metric) on $(\Sigma, T)$ is a metric on $\Sigma$ so that its
restriction to each triangle is isometric to a Euclidean triangle.
It is uniquely determined by a function $\ell : E\rightarrow
\mathbb{R}_{>0}$ which assigns a length to each edge so that the
triangle inequalities hold for every triangle in $F$.
 With a PL metric, the triangulated surface $(\Sigma, T)$ is
locally isometric to the Euclidean plane or half-plane if there is
boundary except at the vertices where the metric may have
cone-like singularities.  In this paper, we will always assume a
triangulated surface is equipped with some PL metric. We denote by
$ij$ the edge with vertices $i$ and $j$, by $ijk$ the triangle
with vertices $i$, $j$ and $k$. If $f$, $g$ and $h$ are the
functions over $V$, $E$, and $F$, respectively, for simplicity, we often
write $f_i$, $g_{ij}$ and $h_{ijk}$ for $f(i)$, $g(ij)$ and
$h(ijk)$.

\begin{definition}[Luo~\cite{luo_yamabe}].
Two PL metrics $\ell$ and $\tilde{\ell}$ on combinatorially
equivalent triangulated surfaces $(\Sigma, T)$  are PL conformal
if
\begin{equation}
\tilde{\ell}_{ij} = e^{w_i + w_j}\ell_{ij}
\end{equation}
 for some function $w: V\rightarrow \mathbb{R}$ and for
all edges $ij$. \label{def:plconformal}
\end{definition}
The function $w$ plays an analogous role of conformal factor in
this PL setting and thus is called a \emph{PL conformal factor}.
We call such change of a PL metric on $(\Sigma, T)$ a \emph{PL
conformal change}. This defines an equivalent relation on PL
metrics on $(\Sigma, T)$. We call an equivalent class \emph{PL
conformal class} on $(\Sigma, T)$.  Motivated by the smooth Yamabe
problem, Luo \cite{luo_yamabe} considered the existence and
uniqueness of PL metrics with prescribed curvature in a PL
conformal class. Namely, suppose $d$ is a PL metric on $(\Sigma,
T)$ and $K: V \to \mathbf R$ is given.  Is there a PL metric $d'$
on $(\Sigma, T)$ which is PL conformal to $d$ so that the
curvature of $d'$ is $K$ (discrete Yamabe problem)?  Is the metric
$d'$ unique up to scaling?  There have been work done on solving
both questions for finite triangulations of compact surfaces.
The main issue that we address in the paper is on a conjecture of 
Luo \cite{luo1} about the uniqueness of the simplest infinite triangulation 
of the plane. Namely,
suppose $d$ is a  PL metric on the hexagonal triangulation of the
plane so that (1) it is PL conformal to the regular hexagonal
tiling, (2) the metric $d$ is complete and has zero curvature at
each vertex. Is $d$ the regular hexagonal tilling?  Our main
result gives an affirmative answer to this question for those PL
metrics so that all inner angles are in $(0, \pi/2-\delta]$ for
some $\delta >0$.

\begin{figure}[t]
\begin{center}
\begin{tabular}{c}
\includegraphics[width=0.4\textwidth]{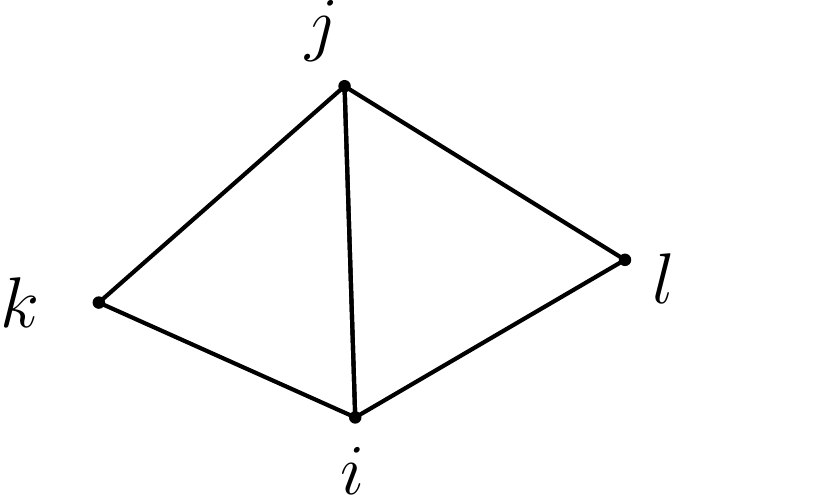}
\end{tabular}
\end{center}
\caption{The length cross ratio of the edge $ij$ is defined as $\frac{l_{il}l_{jk}}{l_{ik}l_{jl}}$.}
\label{fig:cross_ratio}
\end{figure}

One can also look at the PL conformal transformation in terms of
cross-ratio.  See Bobenko, Pinkall and
Springborn~\cite{bobenko_conformal}. For an interior edge $ij$
incident to triangles $ijk$ and $ilj$ as in
Figure~\ref{fig:cross_ratio}, if  the quadrilateral $iljk$ is
embedded in the complex plane
$\hat{\mathbb{C}}=\mathbb{C}\cup\{\infty\}$ and denote the vertex
positions by $z_i, z_j, z_k, z_l$, any conformal map of
$\hat{\mathbb{C}}$ preserves the cross ratio $(z_i, z_j, z_l, z_k)
:= \frac{z_i - z_l}{z_i - z_k}/\frac{z_j - z_l}{z_j-z_k}$. One can
see that the absolute value of this cross ratio is
$\frac{l_{il}l_{jk}}{l_{ik}l_{jl}}$ (called lenght cross ratio by
Bobenko, Pinkall and Springborn~\cite{bobenko_conformal}) which is
preserved by the PL conformal change since scale factors are
cancelled.

Another motivation for Luo's conjecture comes from
Thurston's conjecture on the rigidity of hexagonal circle packing
in the plane. In his famous address \footnote{International
Symposium in Celebration of the Proof of the Bieberbach
Conjecture. Purdue University, March 1985.} ``The Finite Riemann
Mapping Theorem'', Thurston gave a different approach to PL
conformal geometry for triangulated surface using circle packings,
which takes the view of the conformal map preserving infinitesimal
circles. In Defintion~\ref{def:plconformal},  if we let
$\widetilde{l_{ij}}=(e^{w_i}+e^{w_j})$, we obtain a \it circle
packing metric \rm $\widetilde{l_{ij}}$ in the sense of Thurston.
Thurston conjectured, among other things, that the only complete
flat circle packing metric on the hexagonal triangulation of the
plane is the regular hexagonal packing.  This was proved by Rodin
and Sullivan~\cite{RS}. The main problem that we study is the
counterpart of Thurston's conjecture in the new  PL conformal
setting.

In the circle packing setting, it is natural to assign each vertex
$i$ a circle with the radius $e^{w_i}$, and then the circles
centered at two neighboring vertices are tangential to each other.
This PL conformal transformation map circles to circles, or in an
equivalent way. The identity $l_{il}+l_{jk}=l_{ik}+l_{jl}$ always
holds for any edge $ij$. This identity also holds for conformal
transformation on $\hat{\mathbb{C}}$. One disadvantage of this
approach is that in general the meshes we get from practice do not
satisfy $l_{il}+l_{jk}=l_{ik}+l_{jl}$, and usually we have to give
up its edge-length information by assuming $l\equiv const$ to do
the PL conformal transformation. However, in the new PL
conformality the initial metric information can be reserved and
the transformation can be done directly with the original mesh.

Let $\alpha_i$ be the cone angle at the vertex $i$ which is the
sum of all inner angles having the vertex $i$. The curvature is a
function over the vertices: $K: V \rightarrow \mathbb{R}$ defined
by $K_i = 2\pi - \alpha_i$ if the vertex $i$ is in the interior
and $K_i = \pi - \alpha_i$ if the vertex $i$ is on the boundary.
It is obvious that the curvature is uniquely determined by the PL
metric. Its converse is the type of rigidity problem (uniqueness
question) we consider here: is the PL metric in a given \emph{PL
conformal class} uniquely
 determined up to a scaling
(i.e., an Euclidean similar transformation) by the curvature
function? In this paper, we will focus on the triangulated planes
with the curvature function $K \equiv 0$, or the triangulated flat
planes for short. Notice that a triangulated flat plane may not be
 isometric to the Euclidean plane.

\subsection{Problem and main results}

Below is the main problem which was first studied in \cite{luo_yamabe}.


{\it Given two PL metrics $\ell$ and $\tilde{\ell}$ on
combinatorially equivalent triangulated surfaces $(\Sigma, T)$ so
that both metrics are isometric to the complex plane $\mathbb{C}$,
if $\ell$ and $\tilde{\ell}$ are PL conformally equivalent, are
they differ by a scaling?  }



\begin{figure}[t]
\begin{center}
\begin{tabular}{c}
\includegraphics[width=0.5\textwidth]{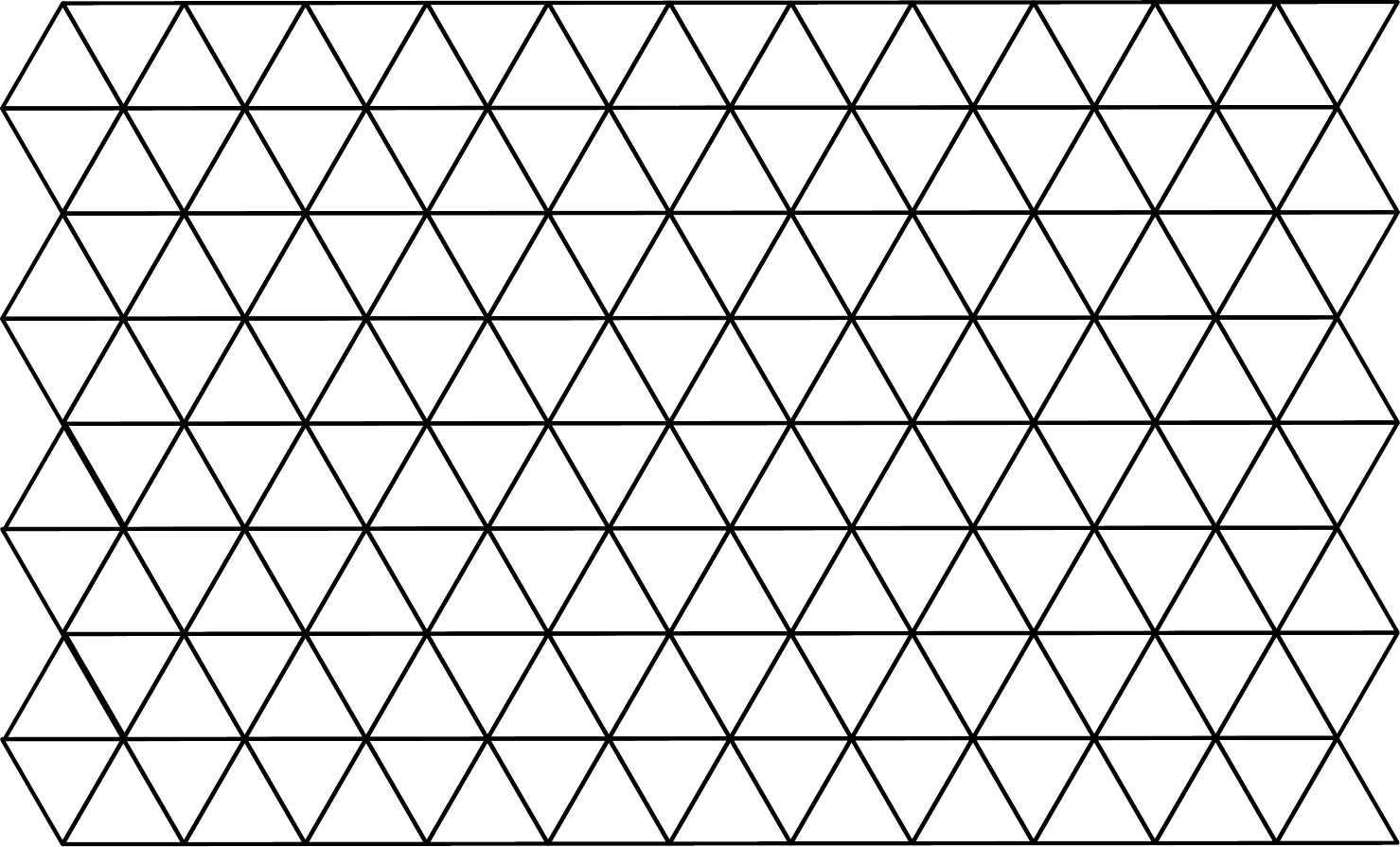}
\end{tabular}
\end{center}
\caption{Regular triangulation of the plane.}
\label{fig:regular}
\end{figure}

In this paper, we restrict ourselves to the hexagonal
triangulation of the plane where every vertex is of degree $6$.
The \it regular hexagonal triangulation \rm is the tiling of the
plane by regular equilateral triangles as shown in
Figure~\ref{fig:regular}, i.e., we assume a regularly triangulated
plane is equipped with a PL metric where all the edges have the
same length. It is obvious that a regularly triangulated plane is
isometric to the complex plane. The main result of this paper is
the following rigidity theorem.
\begin{theorem}
If a hexagonally triangulated plane with a piecewise flat metric $\ell$ satisfies the following conditions
\begin{itemize}
\item[(1)] it is PL conformal to the regularly triangulated plane,
and \item[(2)] it is isometric to the complex plane $\mathbb{C}$,
and \item[(3)] $\sup\{\text{all inner angles of the
triangulation}\} <\pi/2$,
\end{itemize}
then it has to be regular, i.e., $\ell$ is constant, or equivalently, the conformal factor $w$ is constant.
\label{thm:main}
\end{theorem}
Althought our proof of the theorem relies on the condition (3), we
believe that it is not a necessary hypothesis and conjecture the
theorem still holds even without it.

The condition (2) is stronger than that flat condition that
curvature $K=0$. Indeed, there are incomplete flat PL metric on
the plane. The condition (2) is equivalent to the complete
flatness. We call the conformal factor $w$ is linear if it is the
restriction of a linear function on $\mathbb{C}$ to the vertices
in a regular triangulation. One can show that a hexagonally
triangulated plane is flat if it is conformal to a regular one
with a linear PL conformal factor $w$ (see
Lemma~\ref{lemma:flat}).  This gives a two-parameter family of the
flat planes of hexagonal triangulation. Other than the regular
ones, each of them must have an overlap of positive area if it is
immersed into the complex plane.

\begin{figure}[t]
\begin{center}
\begin{tabular}{c}
\includegraphics[width=0.9\textwidth]{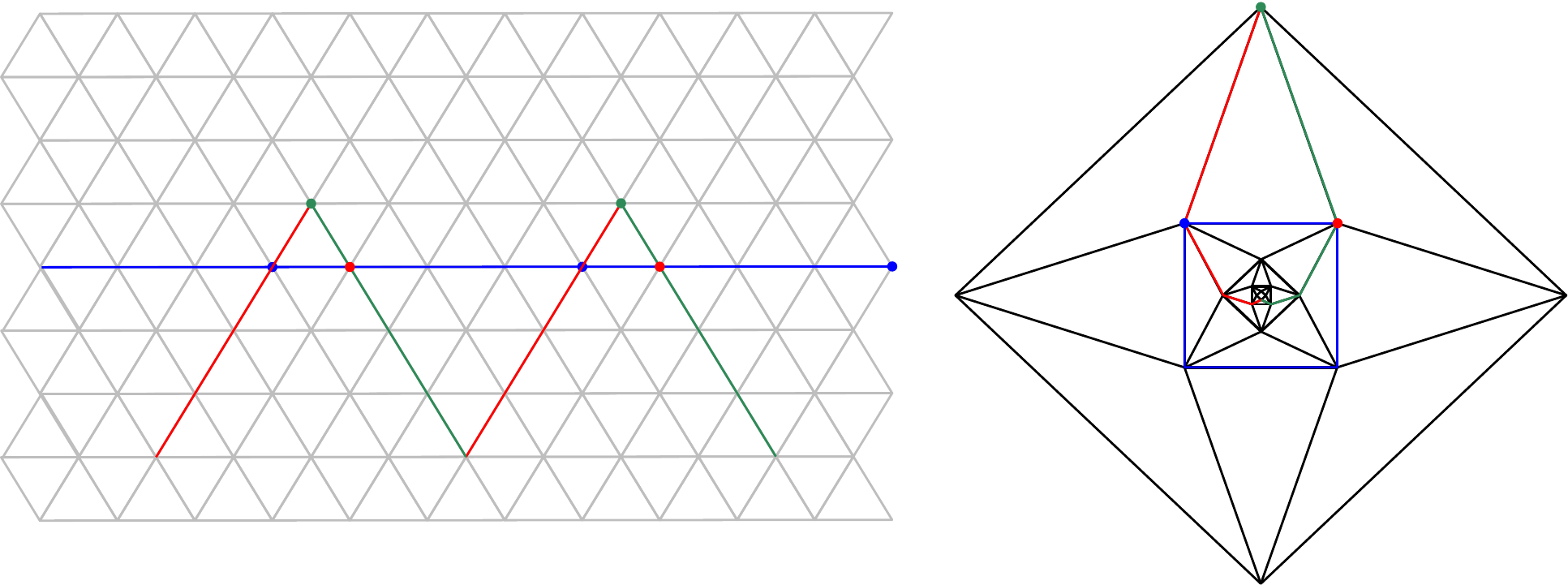}
\end{tabular}
\end{center}
\caption{A hexagonal triangulation of the flat plane.}
\label{fig:hexagonal-flat}
\end{figure}

Figure~\ref{fig:hexagonal-flat}
shows such an example with a linear PL conformal factor $w$. The
left picture shows the regular hexagonal triangulation of
$\mathbb{C}$. The mesh in the right figure contains a nested
infinite sequence of squares so that the $i$th square is
transformed to $(i+1)$th square by a rotation of $\pi/2$ and a
scaling whose ratio is independent of $i$. One can choose this
ratio so that the length cross ratio of each edge in the mesh is
$1$.  This mesh becomes a hexagonal triangulation of
$\mathbb{C}\backslash\{0\}$.  Extend this mesh to the universal
covering space of $\mathbb{C}\backslash\{0\}$, and obtain a flat
plane with the hexagonal triangulation where the length cross
ratio of each edge is $1$. Thus this extension of the mesh in the
right figure is PL conformal to the regularly triangulated plane.
In Figure~\ref{fig:hexagonal-flat}, the colored edges in the left
corresponds to the edges of the same color in the right. One can
verify that the conformal factor $w$ is linear based on the
similarity relation of the triangles in the mesh.  It is shown in
Lemma~\ref{lemma:not_embedded} that this extension of the mesh in
the right figure cannot be immersed to $\mathbb{C}$ without
overlapping. We conjecture that this two parameter family
characterizes all hexagonal triangulation of the flat plane
conformal to the regular one.

\begin{conjecture}
A hexagonally triangulated flat plane is conformal to the regular
one if and only if the PL conformal factor $w$ is linear.
\end{conjecture}

The rigidity of PL conformal transformation for a compact finitely
triangulated surface was initially investigated by
Luo~\cite{luo_yamabe} where he proved that the metric within a PL
conformal class is locally uniquely determined by the curvature by
establishing a variational principle whose action functional is
locally convex function. In 2010,   Bobenko, Pinkall and
Springborn~\cite{bobenko_conformal} found, among other things,  an
explicit formula for the action functional and showed that it
extends to a globally convex one. Using this, they proved a global
rigidity for finite triangulated surfaces. However, the
variational method can not be extended to show the rigidity
problem for infinite triangulated surfaces since action functional
may become infinite.

A similar rigidity results have been proved in the infinite circle
packing case.  Rodin and Sullivan~\cite{RS} and He
~\cite{he_packing} proved the rigidity of hexagonal circle
packings of the complex plane and their methods can be extended to
prove the rigidity of packings with bounded valence.
Schramm~\cite{schramm_rigidity} used a topological property to
extend the rigidity for the packings of arbitrary (locally finite)
planar triangulations. He~\cite{he_rigidity} showed that a
variation along two conformal packings is a harmonic function on a
recurrent network and hence is constant. In this way, he proved
the rigidity of packings even with overlaps.

\section{Outline of the proof}
In this section we outline the proof of Theorem~\ref{thm:main}. With a properly chosen coordinate system,
one can index by a couple of integers $(m, n)$ the vertices in the regular triangulation with the
edge length equals $1$:
\begin{equation*}
V_0=\{m+n\omega|m,n\in \mathbb{Z},\omega=-\frac{1}{2}+\frac{\sqrt{3}}{2}i\}.
\end{equation*}
Since any hexagonal triangulation of the plane $(\Sigma, T)$ has the same combinatorial structure as this
regular one, we also index the vertices $V$ of any hexagonal triangulation in the same manner. We write
$i\sim j$ if $i, j\in V$ are two endpoints of an edge in $E$.

The metric of a hexagonally trangulated plane $(\Sigma,T)$ conformal
to a regular one is uniquely determined by the PL conformal factor $w$
up to a similarity. The condition (2) in Theorem~\ref{thm:main} that
$(\Sigma,T)$ is isometric to $\mathbb{C}$ implies that the curvature
induced by the metric is everywhere $0$, which in return imposes certain
restriction on the function $w$. We introduce a conception called
quasi-harmonic and show the function $w$ is quasi-harmonic.


\begin{definition}
We say a function $f$ on $V$ is quasi-harmonic if there exists $m>0$ depending on $f$ such that,
for any $i\in V$, with its 6 neighbors $i_1,i_2,\cdots,i_6$, there exists $m_1^i,m_2^i,\cdots ,m_6^i\geq m$, satisfying
$$
\sum^6_{j=1}m_j^i=1\quad\quad \text{and}\quad\quad f(i)=\sum^6_{j=1}m_j^if(i_j).
$$
We call such $m$ a harmonic factor of $f$.
\end{definition}

We also define the discrete difference of the PL conformal factor $w$ as follows:
\begin{definition}
For any $c\in V$, we define the difference of $w$ with vector $c$ as $\Delta_cw(i)=w(i+c)-w(i)$.
\end{definition}

With this definition of quasi-harmonic and the difference operator, we have the following lemma:

\begin{lemma}
If a hexagonally triangulated flat plane $(\Sigma,T)$ is PL conformal to a regular one and
$$
\sup\{\text{all inner angles of the triangulation}\}=\theta<\pi/2,
$$
then for any constant $c\in V$, the function $\Delta_cw$ is quasi-harmonic and
its harmonic factor $m(\theta)$ depends only on $\theta$.
\label{lemma:quasiharm}
\end{lemma}

By definition if $m$ is a harmonic factor of a quasi-harmonic function $f$, any
$\widetilde{m}\in(0,m)$ is also a harmonic factor. It should be noticed that a
quasi-harmonicity is weaker than harmonicity related to the graph Laplace operator.
In quasi-harmonicity the weight $m_i^j$ is directed, i.e., $m_i^j$ is not necessarily equal to $m_j^i$.
This means the random walk on the $1$-skeleton of $T$ defined by the weights $m_i^j$ may not
be reversible. Nevertheless, a quasi-harmonic function satisfies the maximal principle.
Notice that it is only Lemma~\ref{lemma:quasiharm} whose proof in this paper requires
the condition (3) that all inner angles of the hexagonal triangulation are strictly acute.


It is well-know that a bounded harmonic function on a recurrent network must be constant ~\cite{he_rigidity}(Lemma 5.5).
We show that a quasi-harmonic function on a network has a similar property of almost
constant over an arbitrarily large region. We assume a graph distance between any two
vertices $i, j$ in $V$, i.e.,
$$d(i,j)=\inf\{t\in N|\exists \text{a path with $t$ edges in $E$ connecting $i$ and $j$}\}$$
We denote $B(i,R)=\{j\in V|d(i,j)\leq R\}$. Notice that in the definition we use $\leq$ instead of $<$ and it's a little different from the continuous case.

\begin{lemma}
Given a quasi-harmonic function $f$ on $V$ with harmonic factor $m$, if there exist $M\in\mathbb{R},R>0,\epsilon>0$ and $i\in V$ so that
$$
f(i)\geq M-\epsilon m^R \quad\text{and}\quad f|_{B(i,R)}\leq M,
$$
we have $f|_{B(i,R)}\geq M-\epsilon$.
\label{lemma:const_onefun}
\end{lemma}

In the case where the given quasi-harmonic function $f$ is bounded, we can
choose $M$ as its least upper bound and for any $R,\epsilon>0$ choose vertex $i$ satisfying $f(i)\geq
M-\epsilon m^R$, and then by
Lemma~\ref{lemma:const_onefun} we have $f|_{B(i,R)}\geq M-\epsilon$. This shows
there is an arbitrarily large region (specified by $R$) where a bounded quasi-harmonic
function is a constant up to an arbitrarily small pertubation
(specified by $\epsilon$). Based on Lemma~\ref{lemma:const_onefun}, we can prove a stronger result which says
there is an arbitrarily large region where any two bounded quasi-harmonic functions are
simultaneously constant up to an arbitrarily small pertubation.
\begin{lemma}
Given two bounded quasi-harmonic functions $f_1,f_2$ on $V$,  assume the least upper bound of $f_1$ is $M$. Then for any $R>0,\epsilon>0$,
there exists $N\in \mathbb{R}$ and $i\in V$ such that
$$
M-\epsilon\leq f_1|_{B(i,R)}\leq M
$$
$$
N-\epsilon \leq f_2|_{B(i,R)}\leq N.
$$
\label{lemma:const_twofun}
\end{lemma}

Given a PL conformal factor $w$ on the vertices $V$, we focus on two functions $\Delta_1w$ and $\Delta_\omega w$.
The following lemma claims both functions are indeed bounded.
\begin{lemma}
If a hexagonally triangulated flat plane $(\Sigma,T)$ is PL conformal to a regular one,
then $\sup_{i\sim j}|w(j)-w(i)|<\infty$.
\label{lemma:bounded}
\end{lemma}

This, together with Lemma~\ref{lemma:quasiharm}, shows both $\Delta_1 w$ and $\Delta_\omega w$ are bounded quasi-harmonic
under the hypotheses of Theorem~\ref{thm:main}. Thus there is a large region where both functions are close constants. To see the consequence of this fact, we first investigate that of both $\Delta_1w$ and $\Delta_\omega w$
being exactly constant, or equivalently, $w$ being linear. The following lemma claims that the linearity of $w$ implies the
flatness of the plane.
\begin{lemma}
If it is conformal to a regular one with a linear PL conformal factor, then a hexagonally triangulated plane $(\Sigma , T)$ has to be flat.
\label{lemma:flat}
\end{lemma}
As we mentioned in the introduction section, this shows that there exists a two parameter family of
the flat plane conformal to a regular one where $w$ is induced by any linear function $ax + by$ on the
complex plane with $a, b\in \mathbb{R}$. However, we show that only if the function $w$ is constant can
the flat plane isometrically be embedded into the complex plane.


\begin{lemma}
Assume $(\Sigma,T)$ is a hexagonally triangulated plane and conformal to a regular one, with PL conformal factor $w$, if there exist $M>0,N\in \mathbb{R}$, so that
$$
\Delta_1w\equiv M\quad\quad\text{and}\quad\quad
\Delta_\omega w\equiv N,
$$
then there exists $R(M,N) > 0$ depending on $M,N$ such that for any $i\in V$, $B(i,R)\subseteq V$ cannot be isometrically immersed
into $\mathbb{C}$ without an overlap of positive area.
\label{lemma:not_embedded}
\end{lemma}

Based on the above result, we are able to show the following lemma if the PL conformal
factor $w$ is nearly linear over a large region.


\begin{lemma}
Assume $(\Sigma,T)$ is a hexagonally triangulated flat plane and conformal to a regular one, with PL conformal factor $w$,
and all inner angles of the triangulation are in a compact set $S \subseteq(0,\pi)$.
For any $M>0$, there exists $\epsilon>0, R >0$ only depending on $M,S$
such that for any $N\in\mathbb{R},i\in V$ with
$$
M-\epsilon\leq \Delta_1w|_{B(i,R)}\leq M \quad\quad\text{and}\quad\quad
N-\epsilon \leq \Delta_\omega w|_{B(i,R)}\leq N,
$$
then $B(i,R)\subseteq V$ cannot be isometrically embedded
into $\mathbb{C}$ without an overlap of positive area.
\label{lemma:not_embedded_1}
\end{lemma}

With the above lemmas, we are ready to prove Theorem~\ref{thm:main}.

\noindent{\bf Proof of Theorem~\ref{thm:main}.~~}
By Lemma~\ref{lemma:quasiharm} and Lemma~\ref{lemma:bounded}, $\Delta_1w$ and
$\Delta_\omega w$ are both bounded quasi-harmonic. If $w$ is not constant, we may assume
$M=\sup\{\Delta_1w\}>0$. By the condition (3), we have all inner angles are in $[\pi-2\theta,\theta]\subseteq(0,\pi)$.
Choose $R>0, \epsilon>0$ depending on $M$ and $\theta$ according to Lemma~\ref{lemma:not_embedded_1}.
Once $R$ and $\epsilon$ are chosen, one can choose $N$ and $i$ according to Lemma~\ref{lemma:const_twofun}
such that  the PL conformal factor $w$ satisfies
$$
M-\epsilon\leq \Delta_1w|_{B(i,R)}\leq M,\quad\quad\text{ and}\quad\quad
N-\epsilon \leq \Delta_\omega w|_{B(i,R)}\leq N.
$$
Now by Lemma~\ref{lemma:not_embedded_1} this implies the part of the plane covering $B(i,R)\subseteq V$
cannot be isometrically embedded into $\mathbb{C}$ without an overlap of positive area.
This reaches an contradiction to that $(\Sigma,T)$ is isometric to $\mathbb{C}$.
Therefore the PL conformal factor $w$ must be constant.


\section{Properties of $w$ (Lemma~\ref{lemma:quasiharm} and~\ref{lemma:bounded})}
\subsection{Proof of Lemma~\ref{lemma:quasiharm}}

\begin{figure}[t]
\begin{center}
\begin{tabular}{c}
\includegraphics[width=0.75\textwidth]{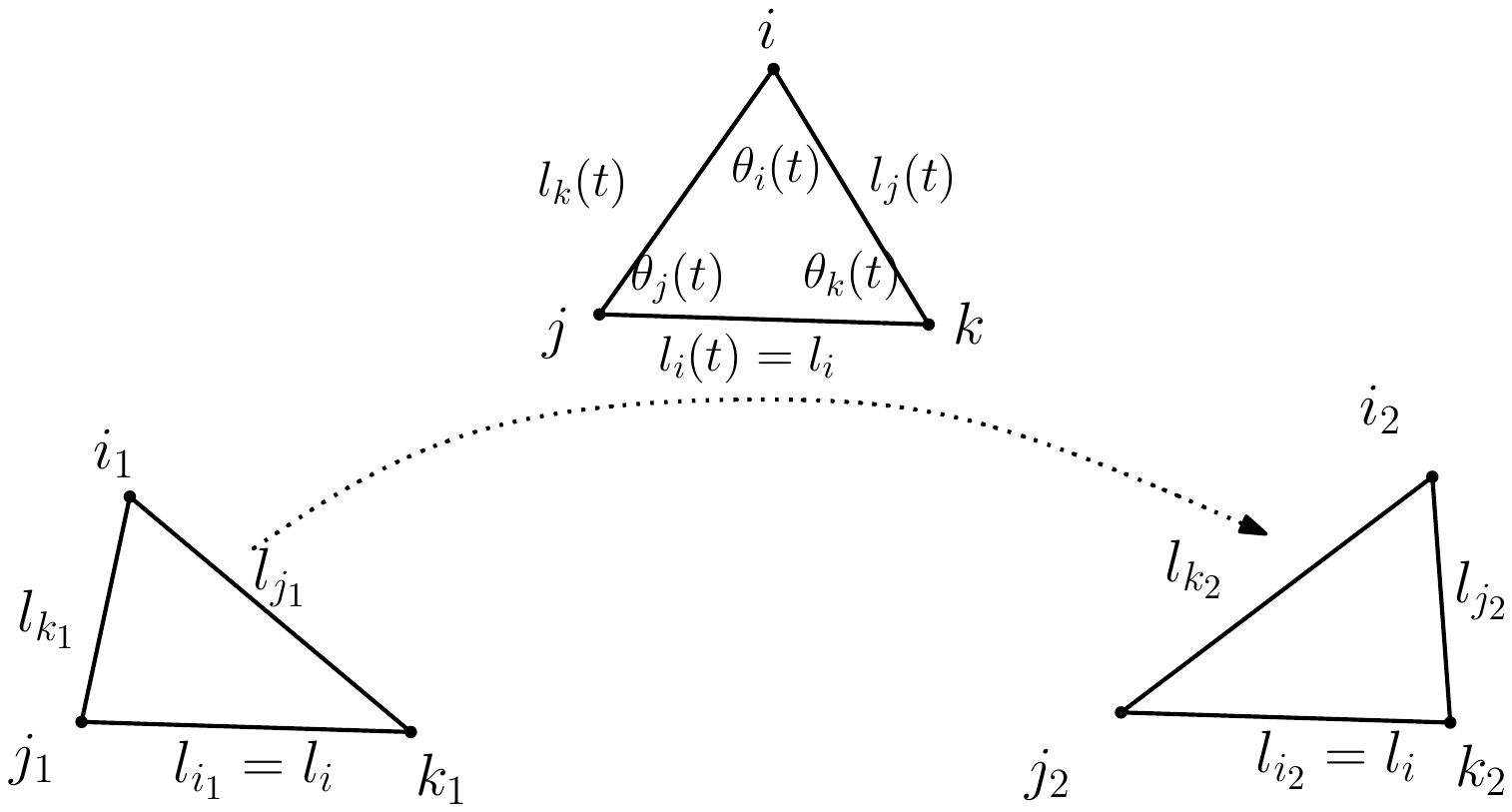}
\end{tabular}
\end{center}
\caption{An acute triangle can be deformed to the target acute triangle in a \textsl{monotonic} way.}
\label{fig:acute_tri_flow}
\end{figure}

We first show the following two lemmas. Lemma~\ref{lemma:acute_tri_flow} shows that when an edge is fixed,
one acute triangle can be deformed to a target acute triangle in a \textsl{monotonic} way. See Figure~\ref{fig:acute_tri_flow}.
Lemma~\ref{lemma:average} describes a relation between the angles and the edge lengths of two triangles.

\begin{lemma}
$\{i_1,j_1,k_1\}$, $\{i_2,j_2,k_2\}$ are two triangles each inner angle
$<\pi/2$ and $l_{i_1}=l_{i_2}$, there exists a flow deforming
$\{i_1,j_1,k_1\}$ to $\{i_2,j_2,k_2\}$ with corresponding edge length $l_i(t)$,
$l_j(t)$, $l_k(t)$ satisfying
\begin{itemize}
\item[(1)] $l_{i}(t)=l_{i_1}=l_{i_2}$;
\item[(2)] $l_{j}(0)=l_{j_1}$, $l_{j}(1)=l_{j_2}$;
\item[(3)] $l_{k}(0)=l_{k_1}$, $l_{k}(1)=l_{k_2}$;
\item[(4)] $l_{j}(t),l_{k}(t), \theta_i(t), \theta_j(t), \theta_k(t)$ are monotonic, continuous and piecewise \\differentiable where
$\theta_i(t), \theta_j(t), \theta_k(t)$ are the 3 inner angles.
\end{itemize}
\label{lemma:acute_tri_flow}
\end{lemma}

\begin{lemma}
For any angle $\theta<\pi/2$ there exist $0<m(\theta)<M(\theta)$ such that for any two triangles $\{i_1,j_1,k_1\}$, $\{i_2,j_2,k_2\}$ with all the inner angle no larger than $\theta$ and edge length $l_{i_2}=l_{i_1}$, $l_{j_2}=l_{j_1}e^{\tilde{u}_j}$, $l_{k_2}=l_{k_1}e^{\tilde{u}_k}$, there exist $a,b$ satisfying $m(\theta)\leq a,b\leq M(\theta)$ and $\theta_{i_2}-\theta_{i_1}=-a\tilde{u}_j-b\tilde{u}_k$.
\label{lemma:average}
\end{lemma}

We will prove Lemma~\ref{lemma:quasiharm} in the following three steps: (1) show Lemma~\ref{lemma:quasiharm} by assuming
Lemma~\ref{lemma:average} holds; (2) show Lemma~\ref{lemma:average} by assuming Lemma~\ref{lemma:acute_tri_flow} holds;
(3) show Lemma~\ref{lemma:acute_tri_flow}.

\subsubsection{Proof of Lemma~\ref{lemma:quasiharm} provided Lemma~\ref{lemma:average} holds}
\begin{figure}[t]
\begin{center}
\begin{tabular}{c}
\includegraphics[width=0.75\textwidth]{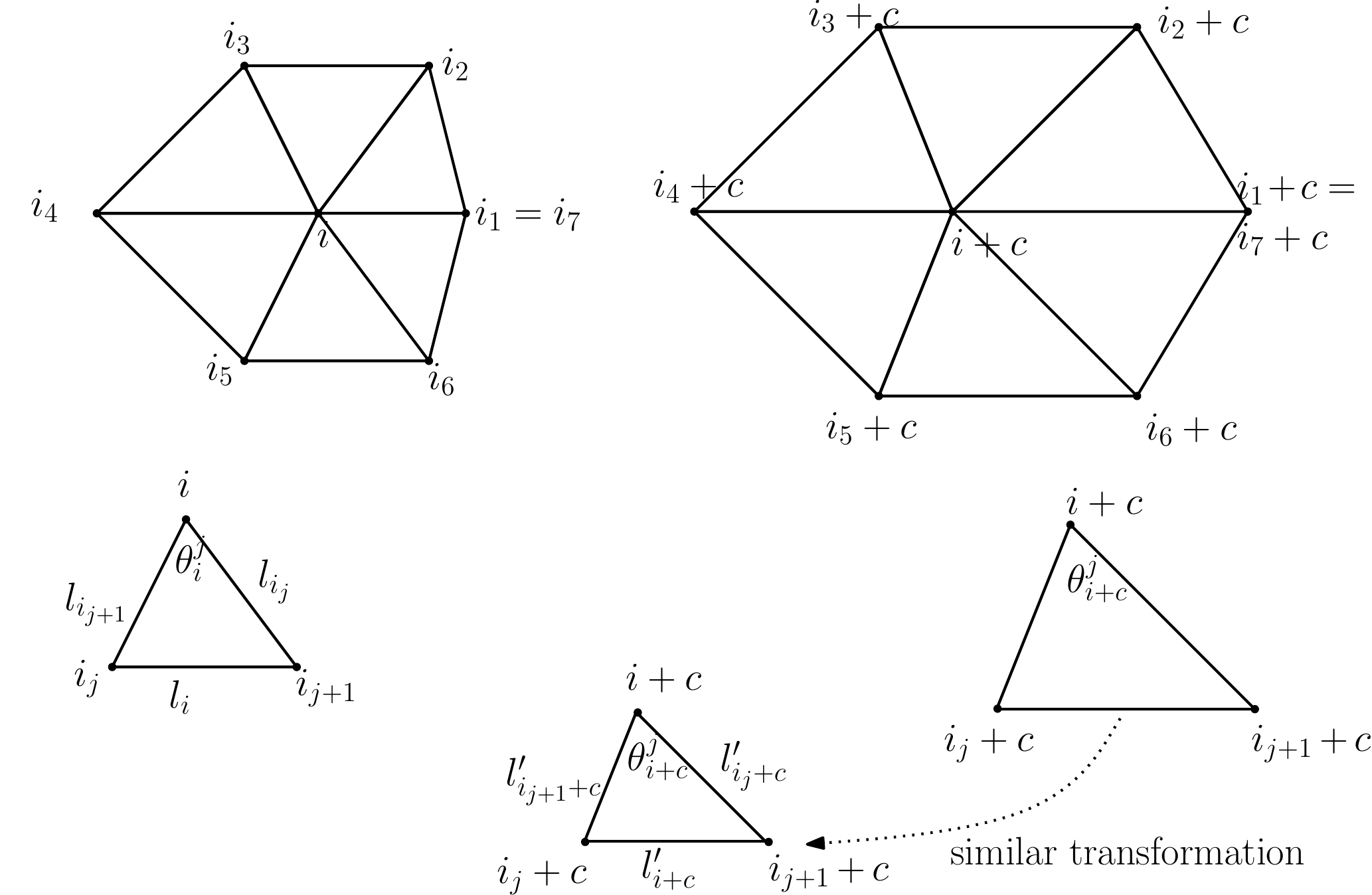}
\end{tabular}
\end{center}
\caption{An acute triangle can be deformed to the target acute triangle in a \textsl{monotonic} way.}
\label{fig:quasi_harmonic}
\end{figure}

Assume $i_1,i_2,\dots,i_6$ are 6 neighbors of $i$ arranged counter-clockwisely, and then $i_1+c,i_2+c,\dots,i_6+c$ are 6 neighbor s of $i+c$ arranged counter-clockwisely. We denote the angle of $i$ in triangle $\{i,j,k\}$ as $\theta_i^{jk}$. For simplicity when $i$ is fixed, we denote $\theta_i^{i_ji_{j+1}}=\theta_i^j$, $\theta_{i+c}^{(i_j+c)(i_{j+1}+c)}=\theta_{i+c}^j$ (assume $i_7=i_1$). See Figure~\ref{fig:quasi_harmonic}

We consider two triangles $\triangle ii_ji_{j+1}$, $\triangle (i+c)(i_j+c)(i_{j+1}+c)$. See Figure~\ref{fig:quasi_harmonic}.  Perform a similar transformation on
the latter so that $l'_{i+c}=l_i$ and  we obtain
\begin{align*}
&l'_{i_{j}+c}=l_{i_{j}} e^{\Delta_c w(i)-\Delta_c w(i_{j+1})}~~\text{and}\\
&l'_{i_{j+1}+c}=l_{i_{j+1}} e^{\Delta_c w(i)-\Delta_c w(i_j)}.
\end{align*}

By applying Lemma~\ref{lemma:average} to $\triangle ii_ji_{j+1}$ and the similarly transformed version of $\triangle (i+c)(i_j+c)(i_{j+1}+c)$,
we have
\begin{align*}
\theta^{j}_{i+c}-\theta^{j}_{i} = -a_j(\Delta_c w(i)-\Delta_c w(i_j))-b_j(\Delta_c w(i)-\Delta_c w(i_{j+1}))\\
\end{align*}
where $m(\theta)\leq a_j,b_j\leq M(\theta)$, $m(\theta),M(\theta)$ are determined by $\theta$ only.
Sum the equality above over all $j=1,2,\dots,6$, we obtain
\begin{align*}
0=&\sum_{j=1}^6(\theta_{i+c}^j-\theta_i^j)\quad\quad\quad\quad\quad\quad\text{(For $K(i)=K(i+c)=0$)}\\
=&-\sum_{j=1}^6(a_j(\Delta_c w(i)-\Delta_c w(i_j))+b_j(\Delta_c w(i)-\Delta_c w(i_{j+1})))
\end{align*}
which leads to
$$
\Delta_c w(i)=\sum_{j=1}^6\frac{a_j+b_{j-1}}{\sum_{j=1}^6(a_j+b_j)}\Delta_c w(i_j))\quad\quad(b_0=b_6)
$$
where
$$\frac{a_j+b_{j-1}}{\sum_{j=1}^6(a_j+b_j)}>\frac{2m(\theta)}{12M(\theta)}=\frac{m(\theta)}{6M(\theta)} \quad\quad \text{and}\quad\quad \sum_{j=1}^6\frac{a_j+b_{j-1}}{\sum_{j=1}^6(a_j+b_j)}=1.$$

Therefore $\Delta_c w$ is quasi-harmonic and $\frac{m(\theta)}{6M(\theta)}$ is its harmonic factor.$\square$

\subsubsection{Proof of Lemma~\ref{lemma:average} provided Lemma~\ref{lemma:acute_tri_flow} holds}
From Lemma~\ref{lemma:acute_tri_flow} there exists a \textsl{monotonic} flow deforming $\{i_1,j_1,k_1\}$ to $\{i_2,j_2,k_2\}$.
Assume $u_j(t)$ and $u_k(t)$ are functions on $[0,1]$ satisfying
$$ l_j(t)=l_{j_1}e^{u_j(t)}\quad\text{and}\quad l_k(t)=l_{k_1}e^{u_k(t)}.$$
It is easy to verify that $u_j(0)=u_k(0)=0$, $u_j(1)=\tilde{u}_j$, $u_k(1)=\tilde{u}_k$. It can be calculated that (e.g.~\cite{})
$$
\frac{\partial \theta_i}{\partial u_j}=-\cot\theta_k\quad\quad \text{and} \quad\quad\frac{\partial \theta_i}{\partial u_k}=-\cot\theta_j.
$$

By the properties of the deforming flow stated in Lemma~\ref{lemma:acute_tri_flow}, we have that
$u_j(t),u_k(t)$ are monotonic, continuous and piecewise differentiable, and
all the inner angles remain in the interval $[\pi-2\theta,\theta]$ during the flow.
Thus we can apply the integral mean value theorem and obtain
\begin{align*}
\theta_{i_2}-\theta_{i_1} =& \theta_i(1)-\theta_i(0)\\
=&\int^1_0\frac{d\theta_i}{dt}dt\\
=&\int^1_0(\frac{\partial \theta_i}{\partial u_j}\frac{du_j}{dt}+\frac{\partial \theta_i}{\partial u_k}\frac{du_k}{dt})dt\\
=&\int^1_0(-\cot\theta_k)\frac{du_j}{dt}dt+\int^1_0(-\cot\theta_j)\frac{du_k}{dt}dt\\
=&-\cot \theta_k(x_k)\int^1_0\frac{du_j}{dt}dt-\cot \theta_j( x_j)\int^1_0\frac{du_k}{dt}dt\quad(0\leq x_j,x_k\leq 1)\quad\\
=&-\cot\theta_k(x_k)u_j(1)-\cot\theta_j(x_j)u_k(1)
\end{align*}
Let  $m(\theta)=\cot\theta$, $M(\theta)=\cot(\pi-2\theta)$.
As $\theta_j(t)$ and $\theta_k(t)$ remain in the interval $[\pi-2\theta, \theta]$ for all $t\in[0,1]$,
we have $m(\theta)\leq\cot\theta_k(x_k),\theta_j(x_j)\leq M(\theta)$. This proves the lemma.$\square$

We remark that here we need the acute triangle assumption to ensure that $\cot\theta_k(t),\cot\theta_k(t)>0$ are positive and bounded.

\subsubsection{Proof of Lemma~\ref{lemma:acute_tri_flow}}
As $l_i$ is fixed in this deforming flow, $\theta_i$, $\theta_j$ and $\theta_k$ are the functions of  $l_j$ and $\l_k$.
Conversely, $l_j$ and $l_k$ are the functions $\theta_j$ and $\theta_k$.
Notice that the triangles are assumed to be acute and therefore $\cot\theta_i,\cot\theta_j$ and $\cot\theta_k$
are all positive. It can be calculated that
\begin{align*}
&\frac{\partial \theta_j}{\partial l_j}=\frac{\cot \theta_i+\cot\theta_k}{l_j}>0, \quad\quad
\frac{\partial \theta_k}{\partial l_j}=-\frac{\cot\theta_i}{l_j}<0, \\
&\frac{\partial \theta_k}{\partial l_k}=\frac{\cot \theta_i+\cot\theta_j}{l_k}>0, \quad\quad
\frac{\partial \theta_i}{\partial l_k}=-\frac{\cot\theta_j}{l_i}<0, \\
&\frac{\partial \theta_i}{\partial l_j}=-\frac{\cot\theta_k}{l_i}<0\quad\quad\text{and} \quad\quad~~~
\frac{\partial \theta_i}{\partial l_k}=-\frac{\cot\theta_j}{l_i}<0.
\end{align*}
Similarly, we have
\begin{align*}
\frac{\partial l_j}{\partial \theta_j}=l_j(\cot \theta_j+\cot\theta_i)>0, \quad\quad
\frac{\partial l_j}{\partial \theta_k}=l_j\cot\theta_i>0, \\
\frac{\partial l_k}{\partial \theta_k}=l_k(\cot \theta_k+\cot\theta_i)>0\quad\text{and}\quad
\frac{\partial l_k}{\partial \theta_j}=l_k\cot\theta_i>0. \\
\end{align*}

We prove the lemma by specifying the deforming flow for different cases as follows.
\begin{itemize}
\item[(1)] If $\theta_j(0)\leq\theta_j(1)$ and $\theta_k(0)\leq\theta_k(1)$, we can choose the flow such that $\theta_j(t)$ and $\theta_k(t)$ are both linear. This means $\theta_i(t)$ is linear too. Thus $\theta'_j(t) \geq 0 $, $ \theta'_k(t)\geq 0$ and $\theta_i(t),\theta_j(t),\theta_k(t)$ are all acute for any $t$.
Therefore $\frac{dl_j}{dt}=\frac{\partial l_j}{\partial \theta_j}{\theta'_j} + \frac{\partial l_k}{\partial \theta_k}{\theta'_k}\geq0$ and $\frac{dl_k}{dt}=\frac{\partial l_k}{\partial \theta_j}{\theta'_j}+\frac{\partial l_k}{\partial \theta_k}{\theta'_k}\geq0$.

\item[(2)] If $\theta_j(0)\geq\theta_j(1)$ and $\theta_k(0)\geq\theta_k(1)$, this case is similar to the case (1).

\item[(3)] If $l_j(1)=l_j(0)$, we choose $l_j(t)=l_j(0)$ and $l_k(t)$ is linear. Now we show that $\theta_i(t),\theta_j(t)$ and $\theta_k(t)$ are monotonous.
It suffices to show they are monotonic in variable $l_k$.
$\frac{\partial \theta_k}{\partial l_k}=(\cot\theta_j+\cot\theta_i)/l_k>0$ as $\theta_i + \theta_j < \pi$. Notice that
we have not shown the acuteness of $\theta_j$ and $\theta_i$. Without loss of generality we may assume that $l_i\geq l_j$,
and thus $\theta_j(t)<\pi/2$. We have $\frac{\partial \theta_i}{\partial l_k}=-\cot\theta_j/l_k<0$ and thus $\theta_i$ is monotonic. This
forces that $\theta_i(t)$ lies between $\theta_i(0)$ and $\theta_i(1)$ and thus  $\cot\theta_i(t)>0$ for any $0<t<1$.
Therefore $\frac{\partial \theta_j}{\partial l_k}=-\cot\theta_i/l_k<0$ and $\theta_j$ is monotonic.

\item[(4)] If $l_k(0)=l_j(1)$, this case is similar to case (3).

\item[(5)] If $\theta_j(0)\leq\theta_j(1)$, $\theta_k(0)\geq\theta_k(1)$, $l_j(0)\leq l_j(1)$ and $l_k(0)\geq l_k(1)$,
we keep angle $\theta_i$ fixed, and let $\theta_j$ increase and $\theta_k$ decrease at the same rate. This can be done
by moving the vertex $i$ along the circumscribing circle of $\triangle i_1j_1k_1$.
Based on the sine law, $l_j$ increases and $l_k$ decreases. All the quantities are moving close to their counterparts
in $\triangle i_2j_2k_2$. Stop once one of following happens: (i)~~$l_j(t)$ reaches $l_j(1)$, or (ii)~~$l_k(t)$ reaches $l_k(1)$, or
(iii)~~$\theta_j(t)$ reaches $\theta_j(1)$, or (iv)~~$\theta_k(t)$ reached $\theta_k(1)$.
Notice that one of the above four cases must happen at some point $t$,  which is a previously discussed case (1 or 2 or 3 or 4).
\item[(6)] If $\theta_j(0)\geq\theta_j(1)$, $\theta_k(0)\leq\theta_k(1)$, $l_j(0)\geq l_j(1)$ and $l_k(0)\leq l_k(1)$, this case is similar to case (5).

Here we remark that the case where $\theta_i(0)=\theta_i(1)$ must falls into either case (5) or case (6).

\item[(7)]If $\theta_j(0)>\theta_j(1)$, $\theta_k(0)<\theta_k(1)$, $l_j(0)< l_j(1)$ and $l_k(0)< l_k(1)$, we fix $l_j$ and increase $\theta_k$. From the cosine law
we have $l_k$ increases at the same time. By the same reasoning as in case (3), both $\theta_i$ and $\theta_j$ decrease at the same time.
We stop deforming once one of the following happens: (i)~~$\theta_i$ reaches $\theta_i(1)$, or (ii)~~$\theta_j$ reaches $\theta_j(1)$,
or (iii)~~$\theta_k$ reaches $\theta_k(1)$, or (iv)~~$l_k$ reach $l_k(1)$. Again, During the deformation,
All the quantities are moving close to their counterparts in $\triangle i_2j_2k_2$. Thus one of the above four cases must happen at some point $t$,
which is a previously discussed case.

\item[(8)] There are three cases left: (i)~~$\theta_j(0)<\theta_j(1)$, $\theta_k(0)>\theta_k(1)$, $l_j(0)> l_j(1)$ and $l_k(0)> l_k(1)$, or
(ii)~~$\theta_j(0)<\theta_j(1)$, $\theta_k(0)>\theta_k(1)$,  $l_j(0)< l_j(1)$ and $l_k(0)< l_k(1)$, or
(iii)~~ $\theta_j(0)>\theta_j(1)$, $\theta_k(0)<\theta_k(1)$, $l_j(0)> l_j(1)$ and $l_k(0)> l_k(1)$.
They are similar to  case (7).$\square$
\end{itemize}

\subsection{Proof of Lemma~\ref{lemma:bounded}}

\begin{figure}[t]
\begin{center}
\begin{tabular}{c}
\includegraphics[width=0.4\textwidth]{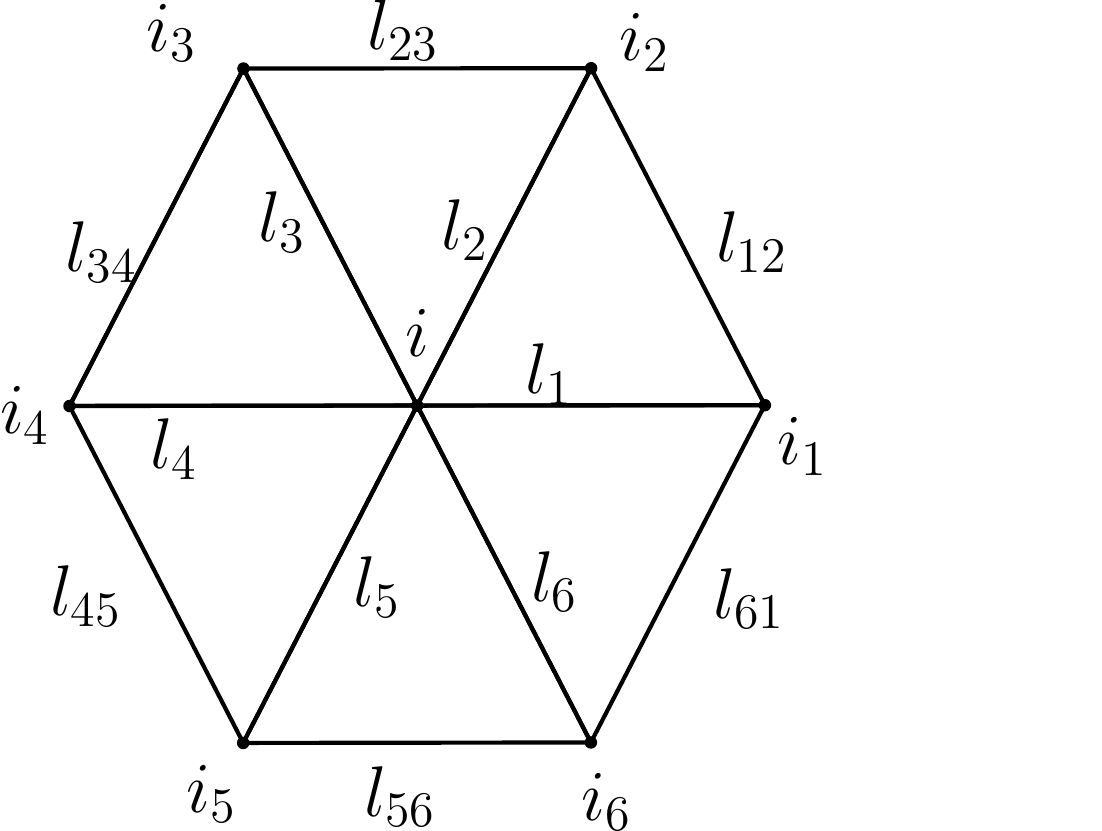}
\end{tabular}
\end{center}
\caption{Illustration for the proof of Lemma~\ref{lemma:bounded}. }
\label{fig:star_1}
\end{figure}

We prove $\Delta_1w(i)=w(i+1)-w(i)$ has a lower bound and the other inequalities can be shown in a similar way.
Let $i_1,i_2,\dots,i_6$ be six neighbors of $i$ and $i_1=i+1$, see Figure~\ref{fig:star_1}.
Let $e^{w(i_1)-w(i_{})} = 1/m$, we show that $m$ is no bigger than $6$.
Consider the triangle $\triangle ii_1i_2$. We have $l_{12}/l_1 = e^{w(i_1)-w(i)}= 1/m $
By triangle inequality, $l_2 > (m-1)l_{12}$, and thus
$e^{w(i_2)-w(i)} = l_{12}/l_{2} < 1/(m-1)$. Now consider the
triangle $\triangle ii_2i_3$, similarly we obtain $e^{w(i_3)-w(i)} < 1/(m-2)$.
We can continue to consider the triangles around $i$ and
obtain
\begin{align*}
e^{w(i_4)-w(i)}<1/(m-3), \dots,~\text{and}~e^{w(i_6)-w(i)}<1/(m-5).
\end{align*}
By contradiction,  assume $m>6$. For any $1\leq j \leq 6$, we have $e^{w(i_j)-w(i)}<1$ and $w(i)>w(i_j)$.
Thus $l_{j,j+1}$ is the shortest side in $\triangle ii_ji_{j+1}$ and
its corresponding angle $\theta_i^j$ is the smallest and thus less than $\pi/3$.
Then $K_i=2\pi-\sum^6_{j=1}\theta_i^j>0$. It contradicts to $K\equiv0$.
Therefore $m\leq6$ and $e^{w(i_1)-w(i)}=1/m>1/6$, and $w(i_1)-w(i)\geq \log(1/6)$. $\square$

\section{Properties of Quasi-harmonic Function (Lemma~\ref{lemma:const_onefun} and~\ref{lemma:const_twofun})}
\subsection{Proof of Lemma~\ref{lemma:const_onefun}}
Consider any $j\in B(i, R- 1)$, let $j_1,j_2,\cdots,j_6$ be its six neighbors.
Since $f$ is quasi-harmonic, there exist $m_1,m_2,\cdots ,m_6\geq m$,
satisfying $\sum^6_{k=1}m_k=1$ and $f(j)=\sum^6_{k=1}m_kf(j_k)$. Thus for any
$k$, we have $M - f(j_k) \leq (M-f(j))/m$. In other words, for any two
neighboring vertex $j\backsim l$ with $j\in B(i, R- 1)$, we have
$M-f(l)\leq(M-f(j))/m$. By induction, we can show that $M-f(j)\leq(M-f(i))/m^n$
if $d(i,j)=n\in \mathbb{N}^+$. In particular, since $M-f(i)<\epsilon m^R$,
for any $j\in B(i,R)$, $M-f(j)<\epsilon$, i.e., $M-\epsilon\leq f|_{B(i,R)}$.
$\square$

\subsection{Proof of Lemma~\ref{lemma:const_twofun}}
Let $m$ be harmonic factor of both $f_1$ and $f_2$. Choose a proper $M_2$ such that $|f_2(i)|<M_2$ for any vertex $i$.
Let $n$ be an integer larger than $2M_2/(\epsilon m^R)$ and $R_2=nR$. Since $M$ is the least upper bound of $f_1$,  there exists a vertex
$i$ such that $f_1(i)>M-\epsilon m^{R_2}$. By Lemma~\ref{lemma:const_onefun}, $f_1|_{B(i,R_2)}>M-\epsilon$.

Let $F(k)$ be the maximum of $f_2$ in $B(i,kR)$, i.e. $F(k)=\max_{j\in B(i,kR)}f_2(j)$ ($k=0,1,\dots,n$). By the definition of $F$, we have $-M_2\leq F(0)\leq F(1)\leq\dots\leq F(n)\leq M_2$.
So there exists $k\in\{1,2\dots,n\}$ such that $F(k)-F(k-1)\leq 2M_2/n\leq\epsilon m^R$. Choose $j_0\in B(i,(k-1)R)$ s.t. $f_2(j_0)=F(k-1)\geq F(k)-\epsilon m^R$. For $B(j_0,R)\subseteq B(i,kR)$, we have $f_2|_{B(j_0,R)}\leq F(k)$. By Lemma~\ref{lemma:const_onefun},
$f_2|_{B(j_0,R)}\geq F(k)-\epsilon$. Let $N=F(k)$ and we have $N-\epsilon \leq f_2|_{B(j_0,R)}\leq N$. As $B(j_0,R)\subseteq B(i,R_2)$, we also have $M-\epsilon\leq f_1|_{B(j_0,R)}<M$. We can just take $j_0$ here as $i$ in the lemma.$\square$

\section{ Consequences of $w$ being (almost) linear (Lemma~\ref{lemma:flat}, \ref{lemma:not_embedded} and~\ref{lemma:not_embedded_1})}
\subsection{Proof of Lemma~\ref{lemma:flat}}
\begin{figure}[t]
\begin{center}
\begin{tabular}{c}
\includegraphics[width=0.7\textwidth]{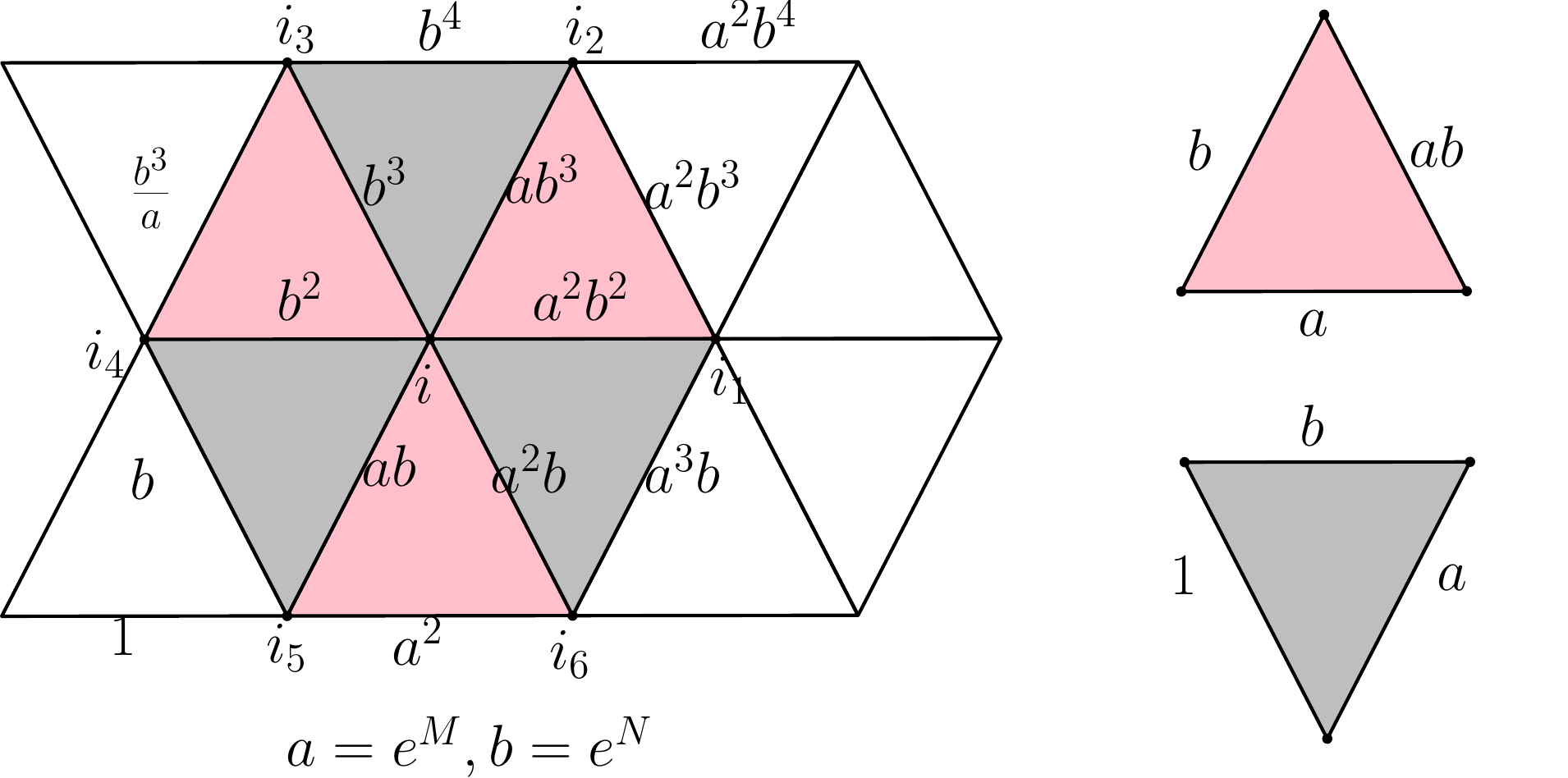}
\end{tabular}
\end{center}
\caption{In the case where $w$ is linear, there are only two types of triangles in the hexagonal mesh up to a similar transformation.}
\label{fig:linearity_flat}
\end{figure}

Since $w$ is linear, there exists $M$  and $N$ such that $\Delta_1 w \equiv M$
and $\Delta_\omega w \equiv N$.  One can show that for any vertex $i$,
$\triangle ii_1i_2$, $\triangle i_4ii_3$ and  $\triangle i_5i_6i$ are
similar to each other,  and thus $\theta^1_i+\theta^3_i+\theta^5_i=\pi$. See
Figure~\ref{fig:linearity_flat}. For the same reason,
$\theta^2_i+\theta^4_i+\theta^6_i=\pi$. Therefore  $K_i=0$ for any vertex $i\in
V$.

\subsection{Proof of Lemma~\ref{lemma:not_embedded}}
\begin{lemma}
Assume we have an immersion $g$ from $\Sigma$ to $\mathbb{C}$, and we take
vertex in $V$ just as the points in $\Sigma$, unit vector $e\in\{\pm1,\pm
\omega,\pm(\omega+1)\}$, for any $i\in V$, we have a unique orientation
preserving similar transformation (non-degenerated 1-dimension complex affine
transformation) $T$ on $\mathbb{C}$ such that for any $k\in\mathbb{Z}$,
$T(g(i+ke))=g(i+(k+1)e)$.
\label{lemma:map}
\end{lemma}

Proof: From the proof of Lemma~\ref{lemma:flat}, there are two types of triangles up to similar transformation in the triangulation,
which means the angle between vector $$\overrightarrow{g(i+(t-1)e)g(i+te)},\overrightarrow{g(i+te)g(i+(t+1)e)}$$
is independent of $t$. Here $t$ is an integer.
By the linearity of $w$, one can also verify that $$|\overrightarrow{g(i+(t-1)e)g(i+te)}|/|\overrightarrow{g(i+te)g(i+(t+1)e)}|$$ is
independent of $t$. Thus there exists the unique $k\in \mathbb{C}$ such that $$(g(i+(t+1)e)-g(i+te))=k(g(i+te)-g(i+(t-1)e)).$$
If denote $b = g(i+(t+1)e)-kg(i+te)=g(i+te)-kg(i+(t-1)e)$, then $b$ is a constant independent of $t$.
Therefore $g(i+(t+1)e)=kg(i+te)+b$, and the transformation $T:z\mapsto kz+b$
is the unique map as we claimed.$\square$

Assume we have an immersion $g$ from $\Sigma$ to $\mathbb{C}$. From the hypothesis, $w$ is linear on $V_0$.
For any $i\in V$, there are at least two unit vectors $e_1,e_2\in\{\pm1,\pm \omega,\pm(\omega+1)\}$
such that $w(i+e_1)-w(i)<0$, $w(i+e_2)-w(i)<0$. By Lemma~\ref{lemma:map}, there exists the unique contract
affine transformation $T_1$ on $\mathbb{C}$ which maps $g(i+te_1)$ to $g(i+(t+1)e_1)$ and
has the unique fixed point denoted $O_1$. Similarly, there exists the unique contract affine
transformation $T_2$ on $\mathbb{C}$ which maps $g(i+te_2)$ to $g(i+(t+1)e_2)$ and
has the unique fixed point denoted $O_2$.

For $||T_1||=e^{2(w(i+e_1)-w(i))}<1,||T_2||=e^{2(w(i+e_2)-w(i))}<1$, by the
fixed point theorem $T^m_1(g(i))\rightarrow O_1$, $T^m_2(g(i))\rightarrow O_2$
as $m\rightarrow\infty$.
\begin{align*}
&|O_1-O_2|\\
=&\lim_{m\rightarrow\infty}|T^m_1(g(i))-T^m_2(g(i))|\\
=&\lim_{m\rightarrow\infty}|g(i+me_1)-g(i+me_2)|\\
\leq&\lim_{m\rightarrow\infty}(\sum_{s=1}^m|g(i+me_1+(s-1)e_2)-g(i+me_1+se_2)|\\
&+\sum_{s=1}^m|g(i+me_2+(s-1)e_1)-g(i+me_2+se_1)|)\\
\leq&\lim_{m\rightarrow\infty}(m|g(i+me_1)-g(i+me_1+e_2)|+m|g(i+me_2)-g(i+me_2+e_1)|)\\
=&\lim_{m\rightarrow\infty}m(||T_1||_2^m|g(i)-g(i+e_2)|+||T_2||_2^m|g(i)+g(i+e_1)|)\\
=&0  \quad\quad\quad\quad(\text{For }||T_1||,||T_2||<1)
\end{align*}

So $O_1=O_2$, we may assume it is the origin and $T_1(z)=r_1e^{2\pi
i\alpha_1}z,T_2(z)=r_2e^{2\pi i\alpha_2}z$. There exists
$m_s\rightarrow+\infty, n_s\rightarrow+\infty$ such that
$r_1^{m_s}r_2^{-n_s}\rightarrow1$. And for any $\epsilon>0$, because
$\{|T_2^{-n_s}T_1^{m_s}(g(i))|:\forall s\}$ is bounded, there exists $s\neq t$
such that $|T_2^{-n_s}T_1^{m_s}(g(i))-T_2^{-n_t}T_1^{m_t}(g(i))|<\epsilon/2$
and $r_1^{m_s}r_2^{-n_s}>1/2$, and thus
$|g(i)-T_2^{n_s-n_t}T_1^{m_t-m_s}(g(i))|<\epsilon$. Choose
$R=|n_s-n_t|+|m_t-m_s|$ and then $T_2^{n_s-n_t}T_1^{m_t-m_s}g(i)\in B(i,R)$,
when $\epsilon$ is small enough, $T_2^{n_s-n_t}T_1^{m_t-m_s}g(i)$ must be in
the hexagonal neighborhood of $g(i)$ and this indicates overlapping.$\square$

\subsection{Proof of Lemma~\ref{lemma:not_embedded_1}}
For the sake of convenience,  assume $i=0$ and  $w(0)=0$. For otherwise one can perform a similar transformation to the mesh to make $w(0)=0$.
For any positive integer $R$, let
\begin{align*}
W_R&=\{w:B(0,R)\rightarrow \mathbb{R}| w(0)=0\}, \\
W_R^f&=\{w\in W_R|\text{$w$ induces zero curvature in $B(0,R)$}\}, \text{and}\\
W_R^o&=\{w\in W_R^f|\text{The immersion of $B(0,R)$ into the plane has an overlap of positive area}\}.
\end{align*}

${W_R}$ is a linear space of finite dimension and can be equipped with a metric, for instance, induced from $L_2$ norm.
$W_R^f$ and $W^o_R$ are two subsets (not necessary subspaces) of $W_R$ and have a natural inherent metric.
Based on Lemma~\ref{lemma:bounded}, $W_R^f$ is bounded. When we immerse the mesh in $B(0, R)$ with $w\in W_R^f$
into the plane, we have the freedom of choosing a base point and the orientation of an edge incident to the base point.
On the other hand, once they are chosen, the immersion is uniquely determined. By induction, one
can show that the positions of the immersed vertices are continuous functions of $w$. In addition,
whether the immersion has an overlap of position area is independent of the choice of the base point
and the orientation of that edge. Therefore for any $w\in W^o_R$, any sufficiently small perturbation of $w$ in $W^f_R$
will not move $w$ out of $W_R^o$. So $W_R^o$ is an open subset of $W_R^f$.

According to the hypothesis of the lemma, we fix a $M > 0$. For any $N$, let $w_{M, N}$ be the linear PL conformal factor
i.e., $\Delta_1 w_{M, N}\equiv M$ and $\Delta_\omega w_{M, N}\equiv N$. By Lemma~\ref{lemma:flat}
$w_{M,N}|_{B(0, R)}\in W^f_R$ for any $R$. For simplicity, when it is clear from the context, we also denote $w_{M,N}$
its restriction to $B(0, R)$. By Lemma~\ref{lemma:not_embedded}, we know for any $N$, there exists $R(N)$ large enough
such that $w_{M,N}\in W^o_{R(N)}$. Since $W^o_{R(N)}$ is open in $W^f_{R(N)}$,
there exists a neighborhood of $w_{M,N}$ in $W^f_R(N)$ which remains in $W^o_{R(N)}$.  In particular, there exists $\epsilon(N)>0$
sufficient small such that for any $N'\in (N -\epsilon(N), N+\epsilon(N))$, any $w \in W^f_{R(N)}$ satisfying
$|\Delta_1 w-M|<\epsilon(N)$ and $|\Delta_\omega w-N'|<\epsilon(N)$ is still in $W_{R(N)}^o$.

So far both $\epsilon(N)$ and $R(N)$ depend on $N$. To obtain $\epsilon$ and $R$ independent
on $N$ as claimed in the lemma, our strategy is to show that all possible $N$
form a compact set. Then based on the above results, we have
an open covering of this compact set using the intervals $(N -\epsilon(N), N+\epsilon(N))$
for any $N$. From the compactness, we have a finite subcover and thus
obtain a uniform $\epsilon$ and $R$ independent on $N$.

Let $\widetilde{S}$ denote an open set such that
$S\subseteq\widetilde{S}\subseteq\overline{\widetilde{S}}\subseteq(0,\pi)$,
and define a set
\begin{align*}
N_M=&\{N|\text{all inner angles in the hexagonal mesh with the conformal
factor $w_{M,N}$ are in $\overline{\tilde{S}}$}\}
\end{align*}
We claim that there exists $\epsilon_0\geq0$ such that
for any $\epsilon\leq\epsilon_0$ and any $R>2$, if $N$ is a real number
so that there exists a conformal factor $w$ satisfying the hypotheses of the lemma, in particular
including that
\begin{itemize}
\item[(i)] the metric induced by $w$ is flat and,
\item[(ii)] all inner angles are in $S$ and,
\item[(iii)] $ M-\epsilon\leq \Delta_1w|_{B(i,R)}\leq M$ and,
\item[(iv)] $N-\epsilon \leq \Delta_\omega w|_{B(i,R)}\leq N$,
\end{itemize}
then $N\in N_M$. This indeed shows $N_M$ is the set of all possible $N$.

\begin{figure}[t]
\begin{center}
\begin{tabular}{c}
\includegraphics[width=0.4\textwidth]{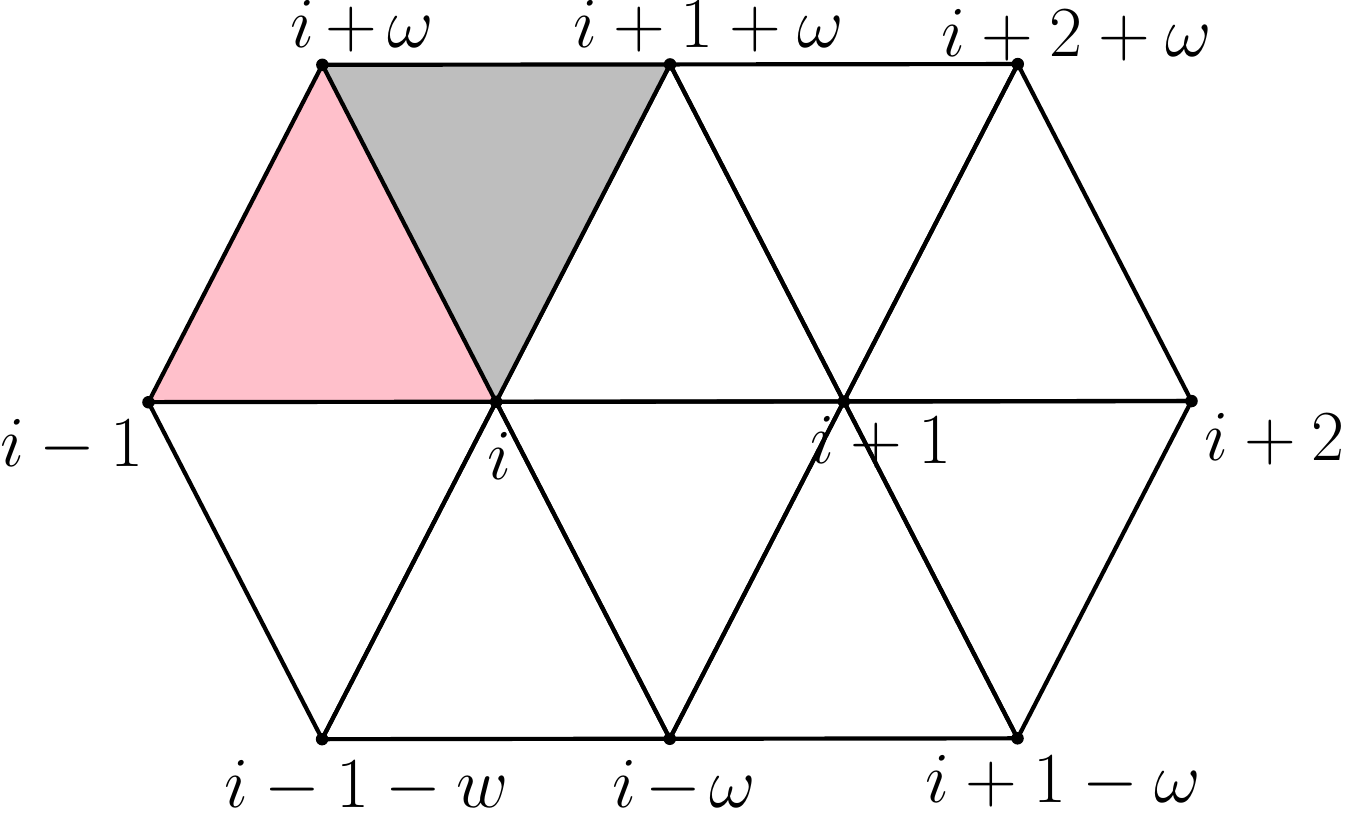}
\end{tabular}
\end{center}
\caption{Illustration for the proof of Lemma~\ref{lemma:not_embedded_1}. }
\label{fig:star}
\end{figure}

To prove the above claim, consider the following two triangles
$\triangle i(i-1)(i+\omega)$ and $\triangle i(i+1+\omega)(i+\omega)$.
See Figure~\ref{fig:star}.  The inner angles
of $\triangle i(i-1)(i+\omega)$ are continuous functions of
$\Delta_1w(i-1)$ and $\Delta_\omega w(i)$, and the inner angles of $\triangle
i(i+1+\omega)(i+\omega)$ are continuous functions of
$\Delta_1w(i+\omega)$ and $\Delta_\omega w(i)$. Denote
\begin{align*}
D_1&=\{(\Delta_1w(i-1),\Delta_\omega w(i))|\text{The inner angles of $\triangle i(i-1)(i+\omega)$ are in } S\} \text{~~and,}\\
D_2&=\{(\Delta_1w(i+\omega),\Delta_\omega w(i))|\text{The inner angles of $\triangle i(i+1+\omega)(i+\omega)$ are all in } S\} \text{~~and,}\\
\widetilde{D_1}&=\{(\Delta_1w(i-1),\Delta_\omega w(i))|\text{The inner angles of $\triangle i(i-1)(i+\omega)$ are all in } \tilde{S}\} \text{~~and,} \\
\widetilde{D_2}&=\{(\Delta_1w(i+\omega),\Delta_\omega w(i))|\text{The inner angles of $\triangle i(i+1+\omega)(i+\omega)$ are all in } \tilde{S}\}.
\end{align*}
By Lemma~\ref{lemma:bounded} we know that $D_1,D_2$ are bounded in $\mathbb{R}^2$.  As $S$ is compact and $\tilde S$ is open,
$D_1$ and $D_2$ are closed and thus compact, and $\widetilde{D_1},\widetilde{D_2}$ are open neighborhoods of $D_1$ and $D_2$ respectively.
So $d(D_1,{\widetilde D_1}^c)>0$, $d(D_2,{\widetilde D_2}^c)>0$.
Choose $\epsilon_0<\frac{1}{2}\min\{d(D_1,{\widetilde {D_1}}^c),d(D_2,{\widetilde {D_2}}^c)\}$.

Now for any $R\geq2$ and any $\epsilon < \epsilon_0$, let $N$ be a number so that
there exists a conformal factor $w$ satisfying the above hypotheses (i, ii, iii, iv).
Then we have $(\Delta_1w(i-1),\Delta_\omega w(i))\in D_1$ and $d((\Delta_1w(i-1),\Delta_\omega w(i)), (M,N)) < 2\epsilon\leq2\epsilon_0 < d(D_1,\widetilde{D_1}^c)$,
Thus $(M,N)\in\widetilde{D_1}$.  Similarly, $(M,N)\in\widetilde{D_2}$.
Notice that from the proof of Lemma~\ref{lemma:flat}, in a hexagonally triangulated plane conformal
to a regular one with a linear conformal factor $w_{N, M}$,
there are only two types of triangles up to similar transformation, which can be represented
by $\triangle i(i-1)(i+\omega)$ and $\triangle i(i+1+\omega)(i+\omega)$. Thus we have $N\in N_M$.

According to the discussion above, for any $N\in N_M$ there exist
$R(N)$ and $\epsilon(N)$ such that for any $N'\in (N -\epsilon(N), N+\epsilon(N))$,
if $w \in W_{R(N)}^f$ satisfies $|\Delta_1w-M|<\epsilon(N)$ and $|\Delta_\omega w-N'|<\epsilon(N)$,
then $w$ is in $W_{R(N)}^o$. Now $\cup_{N\in N_M} (N -\epsilon(N), N+\epsilon(N))$
is an open cover of $N_M$ and there exists a finite subcover $\cup_{j=1}^n B(N_j,\epsilon(N_j))$.
Choose $$\epsilon=\min \{\epsilon_0, \epsilon(N_1), \cdots, \epsilon(N_n) \} \text{~~and~~}
R =\max\{2, R(N_1),\cdots, R(R_n)\}.$$

Finally for any $N$, if $w$ satisfies the above hypotheses (i, ii, iii, iv)
for the chosen $\epsilon$ and $R$,  from the above claim, we have that $N\in N_M$. Thus
$N\in (N_j-\epsilon(N_j), N_j-\epsilon(N_j))$ for some $1\leq j \leq n$. By the choice of $\epsilon$
and $R$, we have $w\in W_{R(N_j)}^f$ and
$|\Delta_1 w-M|<\epsilon(N_j)$ and $|\Delta_\omega w-N|<\epsilon(N_j)$. This implies  $w\in W_{R}^o$ and
proves the lemma. $\square$

\bibliographystyle{plain}
\bibliography{metric_flow}

\end{document}